\documentclass[12pt]{article}
\usepackage{graphicx}
\usepackage{psfrag}
\usepackage{amsmath}
\usepackage{amssymb}
\usepackage{amsthm}

\newtheorem{theorem}{Theorem}
\newtheorem{conjecture}[theorem]{Conjecture}

\newtheorem{definition}[theorem]{Definition}

\newtheorem{hypothesis}[theorem]{Hypothesis}

\begin{document}

\title{Period-doubling cascades for large perturbations of H{\'e}non families}

\author{Evelyn Sander and James A. Yorke}
\date{Preprint \today }

\maketitle

\begin{abstract}
  The H{\'e}non family has been shown to have period-doubling
  cascades.  We show here that the
  same occurs for a much larger class: Large perturbations do not
  destroy cascades. Furthermore, we can classify 
  the period of a cascade in terms of the set of orbits it contains, and count the number of cascades of
  each period. This class of families
  extends a general theory explaining why cascades occur~\cite{sander:yorke:p09}.
\end{abstract}

\section{Introduction}

One of the most mysterious phenomena in nonlinear dynamics is
period-doubling cascades. Cascades never occur alone. Processes have
infinitely many cascades if they have one, and they are seen in a wide
variety of numerical and experimental investigations. As a rule of
thumb, for systems that depend on a parameter, it seems that as
systems become more chaotic, we see period-doubling cascades. We have
taken a more austere view, that if a system is nonchaotic for very
negative values of its parameter and is fully chaotic in some sense
for very large positive parameter values, then for parameters values in
between there are period-doubling cascades, independent of how
complicated the transition is from no chaos to full chaos. By {\bf the
  period} of a periodic orbit, we mean its least period.

A cascade is a special type of connected set of periodic orbits, connected in the space of periodic orbits under the Hausdorff metric. As a first approximation to a definition, we say
a one-parameter family of maps $f:R\times R^n \to R^n$ has a {\bf
  period-$k$ cascade} if there is a (connected) path in $R \times R^n$
of periodic points such that the set of their periods is 
$\left\{ k, 2k, 4k, 8k, \dots \right\}$ (not necessarily occurring in
this order, and listed without multiplicity). 
We restrict these paths to certain ``nonflip'' periodic orbits. The set of periodic orbits can be a collection of complicated networks of periodic orbits, and this restriction prunes the network to manageable simplicity. Specifically, flip orbits are those whose Jacobian matrix has an odd number of real eigenvalues that are less than $-1$ and has none which are equal to $-1$. Nonflip orbits are all the rest. Later, we restate this more formally. 
If the cascade can be chosen so that
its closure in $R\times R^n$ is not compact, then we call the cascade
{\bf unbounded}.  We have developed a theory which explains (and counts) the
occurrence of cascades under general conditions for generic
one-parameter families of $n$-dimensional maps for arbitrary $n$. In
this paper, we show that a family obtained by adding a generic arbitrarily large perturbation to the H{\'e}non family retains the
same cascade structure as the unperturbed family.

The H{\'e}non family \[ H_A(x,y)=\left(
\begin{array}{c}
A + B y -x^2\\
x
\end{array}
\right)
\] 
is a much studied dynamical example.  In an early result in the field,
Devaney and Nitecki~\cite{devaney:nitecki:79} showed that for any
fixed $B$, as $A$ varies from small to large, the H{\'e}non family
forms a horseshoe.  Specifically, all the interesting dynamics is
captured by looking at a certain rectangular region of the plane: For
sufficiently negative $A$ there are no bounded trajectories, but for large positive $A$ the
invariant set in this region has dynamics of a Smale horseshoe. In particular, the
invariant set is a hyperbolic set with one expanding and one contracting
direction, and the dynamics on the set are topologically conjugate to the full shift on two
symbols. Both~\cite{franks:85} and~\cite{yorke:alligood:83} showed
that as the H{\'e}non horseshoe forms, the family has infinitely many
cascades. We now show that this result holds for a broader class of
families.  Some of the difficulties are hinted at in
\cite{yorke:antimonotonicity-Annals92} which shows that the familiar
monotonicity of orbit creation in the one-dimensional quadratic map is
essentially never true for chaotic families of diffeomorphisms in the plane. Orbits are destroyed as well as created as its parameter increases.

We now define the perturbations of the H{\'e}non family 
that we will consider. Note that the word perturbation is usually used to denote something small. In this paper the perturbations can be artitrarily large, as long they are  small in comparison to the original H{\'e}non family in the asymptotic limit as $x,y,$ and $A$ go to infinity. 
We will investigate maps  $F:R \times R^2 \to R^2$ of the form 
\begin{equation}\label{eqn:F} F(A,x,y)= \left(
\begin{array}{c}
A+By-x^2+g(A,x)+\alpha_1(A,x,y)\\
x+\alpha_2(A,x,y)
\end{array}
\right).
\end{equation}
In each case, $B$ is a fixed nonzero constant, and $A$ is the
bifurcation parameter. Functions $g$ and $\alpha=(\alpha_1,\alpha_2)$
satisfy conditions stated below.

The class of functions $g$ permitted are given as follows: Fix $\beta>0$. Define
$G_\beta$ for $\beta > 0$ to be the set of $C^\infty$ functions $g:R^2
\to R$ such that for all $(A,x) \in R^2$,
\[|g(A,0)|<\beta, \mbox{ and } |\partial g/\partial x(A,x)|<\beta. \] 
This class includes for example $C^\infty$ functions that are $C^1$ bounded. 

We now describe the class of functions allowed for $\alpha$. Fix $r>0$
to be any arbitrarily large constant.  For sufficiently small $\delta > 0$ 
depending on $r$, let
\[ \Psi_{\{\delta,r \}} = \{ \alpha: R \times R^2 \to R^2 \in
C^\infty: ||\alpha(A,x,y)||_1< \delta \mbox{ when } ||(A,x,y)|| > r
\} \] (where $||\cdot||_1$ denotes the $C^1$-norm).  Notice that this
class of perturbations has no restrictions other than smoothness in
the region where $||(A,x,y)|| < r$. This class includes for example
all $C^\infty$ functions with compact support, since any $C^{\infty}$
function with compact support is contained in $\Psi_{\{ \delta, r \}}$
for some $r$. Why do we not just assume $\alpha \equiv 0$ outside a
ball or radius $r$? Allowing $\alpha$ to be non-zero everywhere means
that that the set of allowable functions $F$ is open in
$C^\infty$. Thus there exists a residual subset of this open set (the
set depends on $g$) in which all periodic orbit bifurcations are
generic. (See Definition~\ref{def:gen}.) Generic bifurcations allow us
to describe the connected sets of periodic orbits, which is essential
for our task.

Since the function $\alpha$ is uniformly bounded, we can find a constant $\beta_1$ so that for all $(A,x,y) \in R^3$
\begin{equation}\label{eqn:beta} 
|g(A,0)|+|\alpha(A,x,y)|<\beta_1 \mbox{ and } 
|\partial g/\partial x(A,x)|<\beta_1.
\end{equation} 
In a slight abuse of notation, we will drop the subscript $1$, and 
just refer to this new larger constant as $\beta$.

We have chosen these perturbations so that they are dominated by the standard H{\'e}non terms when $A$ and $x$ are large. In particular, for sufficiently large $A = A_1$ and sufficiently small $\delta>0$, the map 
$F(A_1, \cdot, \cdot)$ is topologically the same as the H{\'e}non map. We show that for $A$ equal to  any sufficiently large $A_1$,  any nonflip orbit in the horseshoe  will lie in a cascade, and no other orbit for that $A_1$ is in the cascade. 

\begin{definition}[$PO(f)$ and $PO_{nonflip}(f)$]
  Let $f:R \times R^n \to R^n$ be a $C^\infty$ function. We write
  $[p]$ for the orbit of the periodic point $p$. If $p$ is a periodic
  point for $f(A, \cdot )$, then in a slight abuse of terminology, we
  say that $\sigma=(A, p)$ is a periodic point for $F$.  We write
  $[\sigma]$ or $(A, [p])$ for the orbit.  We denote the set of
  periodic orbits in $R \times R^n$ under the {\bf Hausdorff metric}
  by $PO(f)$.

Let $\sigma=(A, p)$ be a periodic point of period $k$ of a smooth map
$G = f(A,\cdot)$.  We refer to the {\bf eigenvalues of $\sigma$ or
  $[\sigma]$} as shorthand for the eigenvalues of Jacobian matrix
$DG^k(p)$. Of course all the points of an orbit have the same
eigenvalues. We say that $[\sigma]$ is {\bf hyperbolic} if none of its
eigenvalues have norm $1$.

Let $[\sigma]$ be a period-$k$ orbit for $f$.  We call $[\sigma]$ a
{\bf nonflip} orbit for $f$ if $[\sigma]$ has an odd number of real
eigenvalues less than $-1$ and no eigenvalues $= -1$.

We denote the set of all nonflip
  orbits in $R \times R^n$ under the  Hausdorff metric by 
  $PO_{nonflip}(f)$. 
\end{definition}

\begin{definition}[Open arc]
An {\bf open arc} is a set which is homeomorphic to an open interval.
\end{definition}

\begin{definition}[Period-doubling cascade of period $m$]
A  {\bf (period-doubling) cascade of period $m$} is an open arc in 
$PO_{nonflip}(f)$ with the following  properties: 

(i) The open arc contains orbits of
  period $2^k m$ for some positive integer $m$ and for every
  non-negative integer $k$. 

(ii) The number $m$ is the smallest integer for which this is true, and $m$ cannot be made smaller by making the open arc larger.  

Let $(p_k)$ be the sequence of periods of the non-hyperbolic orbits, ordered so that for each $k$, the $k+1$ orbit lies along the open arc between the $k$ orbit and the $k+2$ orbit. Under our genericity hypotheses, it  turns out that no period can occur in the sequence infinitely many times. It follows that in at least one direction ($k \to \infty$ or $k \to -\infty$), $\lim_{k} p_k = \infty$.

We say the cascade is {\bf unbounded} if it does not lie in a compact set of $R\times R^n$.

\end{definition}   

{\bf Even orbits for the two-shift.} For any fixed $k$, consider a
period-$k$ orbit $\overline{S}$ of the full shift on two
symbols. This orbit is associated with a length-$k$ sequence of two
symbols: 
$S = (a_1, \dots, a_k )$, 
where each $a_k$ is
equal to either the symbol $-1$ or $+1$, and $S$ is not periodic.  We say that
$\overline{S}$ is {\bf even} if the associated finite sequence $S$ has
an even number of $-1$'s (or more compactly, if $\Pi_{j=1}^k a_k =1$).

For a map $F(A, \cdot, \cdot)$ of the form in Equation~\ref{eqn:F}, define $MaxInv(A)$ to be the union of trajectories such that  all positive and negative iterates are bounded.

Our main theorem is as follows:
 
\begin{theorem}[Cascades for large perturbations of H{\'e}non families]
\label{maintheorem}
Fix $B \ne 0$, $\beta>0$, and $r>0$. Let $g \in G_\beta$.  For  $\delta>0$, let
$\alpha \in \Psi_{\{ \delta,r \}}$, and let $F$ be as in
Equation~\ref{eqn:F}. 

Then as long as $\delta$ is sufficiently small (depending on $r$), for
every sufficiently large $A=A_1$ depending on $\beta, r$, and $B$, 
there is a residual set of $\alpha
\in \Psi_{\{ \delta,r \}}$ depending on the function $g$ and the constant $B$
for which the following hold:

\begin{enumerate}
\item $MaxInv(A_1)$ is conjugate under a homeomorphism to a two-shift, and this homeomorphism gives a
  one-to-one correspondence of the even symbol sequences with the
  nonflip orbits. (Hence for   $A_1$, we can without confusion refer to a periodic orbit for $F(A_1,\cdot)$ as being even). 

\item Each unbounded cascade contains exactly one periodic orbit for
   $F(A_1,\cdot)$, and it is an even orbit.

\item  For each even orbit there is a unique unbounded cascade containing that orbit. 

\item If an even periodic orbit is period $k$, and $k$ is odd, the
  cascade containing it is a period-$k$ cascade. If $k$ is even, then the cascade containing it is a
  period-$j$ cascade, where $k/j=2^m$ for some $m$.

\end{enumerate}
\end{theorem}

In other words, corresponding to every even period-$k$ symbol
sequence $S$, there is exactly one unbounded cascade of $F$ that satisfies the following:
At parameter value $A_1$, the cascade contains a unique periodic orbit, and it is the unique period-$k$ orbit of $F$ with the symbol sequence $S$.  

{\bf The number of even period-$k$ orbits.} Note that the number of  period-$k$ points of the two shift for each $k$ for has been studied quite extensively, and is often
referred to as the $\zeta$ function. In some cases, it is possible to
write an easy formula for the number $\Gamma (2,k)$ of {\it even} period-$k$
orbits, as follows: There is one even fixed point, and no even
period-2 orbits.  If the $k$ is an odd prime, the number of even
period-$k$ orbits $(2^k-2)/(2k)$. In general, if $k$ is odd, the
number of even period-$k$ orbits is exactly half the number of
period-$k$ orbits. 

In the general case, any positive integer $k$, let  $L(k) = \Sigma( \Gamma(2,j))$ for all $j < k$ for which $k/j$ is a power of $2$. Of course $L(k) = 0$ if $k$ is odd. Then 
\[ \Gamma(2,k) = (\zeta(2,k)/k - L(k))/2. \] 
See~\cite{sander:yorke:p09} for a detailed discussion.

\section{Proof of Theorem~\ref{maintheorem}}

\begin{proof}
Let $B$, $F$, $\alpha$ and $g$ be as in the statement of the theorem. Note that the
assumptions on $g$ and $\alpha$ imply that for all $(A,x,y) \in R^3$, 
\[|g(A,x)|+ |\alpha(A,x,y)| < \beta (1+|x|).\]

Define $s= s(A_1) = \sqrt{A_1}$, and
let $Q=2s$. Let the square $E$ be defined by $E=[-Q,Q] \times [-Q,Q]$.
Assume $A_1 > r$ and $A_1 > Q$. Additional lower bounds will be placed on $A_1$.

The proof of the theorem proceeds from the following steps:

\begin{enumerate}

  \item[\bf Step 1] {\bf Horseshoe dynamics for large $A_1$.} 
Set $\alpha(A,x,y) \equiv 0$.
For  $A_1$ sufficiently large and for $\delta>0$ sufficiently small,  the following are true:

  \begin{enumerate}

  \item[(1a)]  {\bf Periodic orbits in $E$.} For all $A<A_1$, all periodic orbits of $F$ are contained in the
    interior of $E$.

  \item[(1b)] {\bf The two shift.} On $MaxInv(A_1)$, $F$ is topologically conjugate to the full shift on two
    symbols. 
     
  \item[(1c)] {\bf Hyperbolicity.} $F$ is hyperbolic on $MaxInv(A_1)$, with one
    expanding and one contracting direction at each point.
 
   \item[(1d)] {\bf Nonflip orbits.} For $A=A_1$, the nonflip period-$k$ orbits of $F$ are in one-to-one
       correspondence with the even period-$k$ orbits of the full
       shift on two symbols.

  \end{enumerate}

  \item[\bf Step 2] {\bf Adding small perturbations.}
The results in Step 1 are not sensitive to  $C^1$-small perturbations. Thus, they are still true when we add 
$\alpha(A,x,y) \in \Psi_{\{\delta,r \} }$ for sufficiently small $\delta>0$, since we assume that $A_1>r$, implying that $||\alpha||_1 < \delta$ for $A=A_1$. 

\item[\bf Step 3] {\bf No orbits for small $A_0$.} For fixed $A_0$ sufficiently negative (and in particular $A_0 < -r$), the map $F$ has no periodic orbits.

  \item[\bf Step 4] {\bf Cascades.}
Let $\alpha(A,x,y)$ be contained in a residual set of $\Psi_{\{\delta,r\}}$ such that all bifurcations of $F$ are generic (generic bifurcations are defined below). 
Each nonflip periodic orbit of $F(A_1,\cdot,\cdot)$ is contained in a unique unbounded cascade.
   
   \item[\bf Step 5] {\bf Period of the cascades.} 
If $k$ is an odd number, this unbounded cascade is a period-$k$
      cascade of $F$.  If $k$ is an even number, then this unbounded
      cascade is a period-$j$ cascade of $F$, where the ratio $k/j$ is
      a power of two.
\end{enumerate}

\bigskip
{\noindent \bf Proof of Step 1: Horseshoe dynamics for large $A_1$.}

\bigskip
\noindent
{\bf (1a) Periodic orbits in $E$ and (1b) The two shift.}

Let $F$ be of the form in Equation~\ref{eqn:F} with $\alpha \equiv 0$. Let $L=[-Q,Q]$. Let $J_1=[-2s,-s/2]$, $J_2=[s/2,2s]$, $J=J_1 \cup J_2$. 

We have previously shown in~\cite{sander:yorke:p09} that for all sufficiently large parameter values 
$\lambda_1$,
the quadratic map $q(\lambda, x)=\lambda-x^2 + h(\lambda, x)$ -- where
$|h|<\beta (1+|x|)$ -- has the following properties:
\begin{enumerate}
\item[Q1.] $q(\lambda_1,L \setminus J)$ contains no points of $L$.
\item[Q2.] There is an interval $M$ in $L$ such that for all $\lambda \le \lambda_1$,
  each periodic orbit is contained in $M$.
\item[Q3.] At $\lambda_1$, $q(\lambda_1,J_i)$ maps diffeomorphically across
  $L$, where $i=1$ or $2$.
\end{enumerate}

For sufficiently large $A_1$, we get similar results for $F = (F_1, F_2)$:

\begin{enumerate}
\item[F1.] For sufficiently large $A_1$, $F(A_1,E \setminus \{J \times
  L\})$ contains no points of $E$. This is an immediate consequence of
  Q1 above, since $F_1(A_1,x,y)<A_1+|B|Q-x^2$ inside $E$.

\item[F2.] For all $A<A_1$, all periodic orbits are contained in the interior of $E$. This is not immediate from the quadratic case. The proof is as follows:

\begin{proof}
  Let $\{(x_1,y_1) , \dots, (x_k, y_k) \}$ be a periodic orbit at
  parameter $A$.  Fix $x$ to be the $x_i$ with the maximum absolute
  value. Let $y$ be the corresponding $y_i$. Thus $y_i=x_{i-1}$. Let
  $\overline{x}=x_{i+1}$. That is, $F_1(A,x,y)=\overline{x}$. Thus
  $|y|<|x|$, and $|\overline{x}|<|x|$. This implies that,
  \[-|x|<F_1(A,x,y) \le A - x^2 + |B| |y| + \beta(1+|x|).\] 
Since
  $|y|<|x|$, $0 \le (A+\beta) - x^2 + |x| (|B|+ \beta +1)$. Let $\rho
  = (|B|+\beta+1)/2$.  Then 
\[0 \le (A + \beta) + \rho^2 -
  (|x|-\rho)^2.\]
 Hence 
\[|x| \le \rho + \sqrt{A+ \beta + \rho^2}.\]
Note
  that this right-hand side is monotonically increasing in $A$.  Since
  $B$, $\beta$, and thus $\rho$ are fixed, for $A_1$ sufficiently
  large, 
\[\rho+\sqrt{A_1+\beta+\rho^2}<2 \sqrt{A_1}=Q.\] Thus as long
  as $A \le A_1$, we have $|x|<Q$. Since $x$ is the point of the orbit with
  the maximum absolute value, this implies that the periodic orbit is
  contained in the interior of $E$.
\end{proof}

\item[F3.] At $A=A_1$, for each fixed $y$, $F_1(A_1,J_i\times \{ y \} )$
  maps diffeomorphically across $L$. This is immediate from Q3
  above. In addition, $F_2(A_1, J_i \times L)= J_i$. Therefore, each
  $J_i \times L$ maps diffeomorphically, across $E$, and vertically
  staying inside of $E$, such that the $i=1,2$ images are disjoint.
\end{enumerate}

\bigskip
\noindent
{\bf (1c) Hyperbolicity.}

Assume \[\sqrt{A_1}> \beta+|B|+\max \{ 1,|B| \}.\]
The determinant and trace of the Jacobian matrix of $F$ are respectively $-B$
and $-2x + \partial g /\partial x$. For any point in $J\times L$,
\[|-2x + \partial g /\partial x|> 2 s - \beta = 2 \sqrt{A_1} -
\beta.\]
The assumption above on $A_1$  implies that one eigenvalue
for the Jacobian matrix is contracting and the other expanding.

Define the {\it stable} and
{\it unstable cones} respectively by
\begin{eqnarray}
S^+_c &=& \left\{ (\xi,\eta): |\xi| \ge c |\eta| \right\}\\
S^-_c &=& \left\{ (\xi,\eta): |\xi| \le c |\eta| \right\}.
\end{eqnarray}

Then for $A=A_1$ and any point in $J \times L$, the Jacobian matrix
$DF$ maps $S^+_1$ into $S^+_N$, where  $N=\sqrt{A_1}-\beta-|B|$, 
and $DF^{-1}$ maps $S^-_1$ into
$S^-_{N_1}$, where $N_1=N/|B|$. To see this, let
$(\xi,\eta) \in S^+_1$, and let $(\xi_1,\eta_1)=DF(\xi,\eta)$. Then 
\[(\xi_1, \eta_1) = ((-2x+g_x) \xi + B \eta, \xi).\] 
Thus 
\begin{eqnarray*}
\frac{|\xi_1|}{|\eta_1|} & \ge & \frac{ (\sqrt{A_1} - \beta) |\xi| - |B| |\eta| }{|\xi|} \\
& \ge & (\sqrt{A_1}-\beta-|B|)=N>1.
\end{eqnarray*}
Therefore $DF$ maps $S^+_1$ into the interior of itself. 
Likewise, let $(\xi,\eta) \in S^-_1$, and let 
$(\xi_{-1},\eta_{-1})=DF^{-1}(\xi,\eta)$. Then 
\begin{eqnarray*}
\frac{|\eta_{-1}|}{|\xi_{-1}|} &=& \frac{1}{|B|} \frac{|B \xi +(2x - g_x)\eta|}{|\eta|}\\
&\ge& \frac{(\sqrt{A_1}-\beta)|\eta|-|B||\xi|}{|B||\xi|}\\
&\ge& \frac{\sqrt{A_1}-\beta-|B|}{|B|}=N_1>1.
\end{eqnarray*}
Therefore $DF^{-1}$ maps $S^-_1$ into the interior of itself.

Thus  the stable and
unstable cones are mapped strictly inside themselves and expanded
respectively under the derivative and its inverse. 
Using the method of cones  (Corollary 6.4.8
in~\cite{katok:hasselblatt:95}), this guarantees that at $A=A_1$, $F$
is hyperbolic on $MaxInv(F)$.

Putting this together, we get that at $A_1$, $F$ on $MaxInv(F)$ is hyperbolic and topologically conjugate to the two shift. Specifically, we know that $MaxInv(F)$ is contained in $\{ J_1 \cup J_2\} \times L$. The conjugacy codes a point by considering its bi-infinite orbit. For any integer $i$, we code the $i^{th}$ point in the itinerary of an orbit with
a ``$1$'' if the $i^{th}$ iterate is in the left region $J_1 \times L$, and with a ``$-1$'' if the $i^{th}$ iterate is in in the right region $J_2
\times L$.

\bigskip
\noindent
{\bf (1d) Nonflip orbits.}
We now determine the nonflip orbits for $F$ for $A=A_1$. In the left region
$J_1 \times L$, $DF$ has an expanding eigenvalue which is greater than $1$, 
whereas for the right region $J_2 \times L$, $DF$ has an expanding derivative
which is  less than $-1$. Thus any period-$k$ orbit $[p]$ in $J
\times L$ is a nonflip orbit exactly when $[p]$ is in $J_2 \times L$ 
an even number of times. Thus by our conjugacy,
there is the one-to-one correspondence between nonflip
orbits  of $F$ on $MaxInv(F)$ and the
even orbits for the two shift.

\bigskip
{\noindent \bf Proof of Step 2: Adding small perturbations.}

For large $|A_1|$, we have established hyperbolic dynamics on an invariant horseshoe in $MaxInv(F)$. All of this is
 robust under sufficiently $C^1$-small additive perturbations $\alpha(A_1,x,y)$, since $C^1$ small implies that we can make both the function values and all the partial derivatives as small as we want. 
Furthermore, as long as  $\alpha(A,x,y) \in \Psi_{\{\delta,r\}}$, for any $A  \in [A_0,A_1]$, $\alpha$ is $C^1$ small for $|(x,y)|>r>Q$. Therefore $MaxInv(F) \subset E$ for all $A \in [A_0,A_1]$.

\bigskip
{\noindent \bf Proof  of Step 3: No orbits for small $A_0$.}

From F2, it suffices to show that for sufficiently negative
$A_0$, $F(A_0, E) \cap E$ is empty. Note that for all $(x,y) \in E$,
\begin{eqnarray*}
F_1(A_0,x,y)&=&A_0+By-x^2+g(A_0,x,y)+\alpha_1(A_0,x,y)\\
 &<& A_0+\beta + |B| Q  + \beta |x| - x^2 
\end{eqnarray*}
This quadratic in $|x|$ has a maximum at $|x|=\beta/2$, implying that
\[F_1(A_0,x,y) < A_0+\beta+|B|Q + \frac{\beta^2}{4}. \]
Thus as long as $A_0+\beta+(|B|+1)Q +
\frac{\beta^2}{4}<0$, $F_1(A_0,x,y)<-Q$ for any $(x,y) \in E$, implying 
that $F(A_0,x,y)$ is not contained in $E$. 

\bigskip
{\noindent \bf Proof of Step 4: Cascades.}

In order to prove our theorem, we state the following abstract results
on the existence of cascades, from~\cite{sander:yorke:p09}:

\begin{definition}[Generic bifurcations]\label{def:gen}
  Let $f:R\times R^n \to R^n$ be $C^\infty$.  Let $U$ be an open
  subset of $R^{n+1}=R \times R^n$, and let $V$ be its closure.  By
  {\bf periodic orbit bifurcation}, we mean a change (as a parameter 
is varied)  in the local
  number periodic orbits or a change in the dimension of their unstable
  space. We refer to a periodic orbit bifurcation in $U$ as {\bf
    generic} if it is one of the following three types:
\begin{enumerate}
\item A generic saddle-node bifurcation.
\item A generic period-doubling bifurcation.
\item A generic Hopf bifurcation with no eigenvalues which are roots of unity.
\end{enumerate}
\end{definition}

In~\cite{sander:yorke:p09}, we show that a residual set of
one-parameter families have only generic periodic orbit bifurcations.  Let
$\alpha (A,x,y) \in \Psi_{\{\delta, r \}}$ be such that our $F$  has only generic
bifurcations.

We now define the periodic orbit index in a way that is specific to 
$F:R \times R^2 \to R^2$. Let $P$ be a hyperbolic period-$k$ orbit for $F$ with the eigenvalues
$\sigma_1 \le \sigma_2$ for derivative $D(F^k)$ with respect to the
spatial variables $(x,y)$.  Let $I_m=(-\infty,-1)$, $I_0=(-1,1)$, and
$I_p=(1,\infty)$. We define the {\bf periodic orbit index} for $P$ to
be:
\[
ind_F(P) = 
\begin{cases}
1& \mbox{ if } \sigma_1,\sigma_2 \mbox{ are both in } I_m,I_0, \mbox{ or } I_p, \mbox{ or if } \sigma_1,\sigma_2 \mbox{ complex}, \\
-1& \mbox{ if } \sigma_1 \in I_0, \sigma_2 \in I_p, \\
0& \mbox{ if }  \sigma_1\in I_m, \sigma_2 \notin I_m. 
\end{cases}
\]

Note that a flip orbit corresponds to the case of index $0$. For
large parameters such as $A_1$, all periodic orbits are saddles,
implying that the periodic orbit index is -1 or 0.

This definition generalizes to a general definition for  $f: R \times R^n \to R^n$. 
It is a topological invariant, as is described in more generality
and detail in~\cite{sander:yorke:p09}. Let $[\Gamma]:(0,1) \to
PO_{nonflip}(f)$ map homeomorphically to an open arc $C$ in the nonflip 
orbits for $f$. (The brackets are to emphasize that $[\Gamma]$ maps a point in
the interval $(0,1)$ to an orbit in the set of nonflip orbits.) Then
$[\Gamma]$ can be identified with one of the two orientations on $C$.
There is one orientation that is induced by the periodic orbit index:
$[\Gamma]$ is an {\bf index orientation} on $C$ as long as it has the
following property for every $s \in (0,1)$: $ind_f([\Gamma(s)])=-1$
whenever the parameter $A$ is locally decreasing, and
$ind_f([\Gamma(s)])=+1$ whenever the parameter $A$ is locally
increasing. As long as  $f$ has only generic
bifurcations, there exists an index orientation on every open arc in the
set of nonflip periodic orbits.

\begin{hypothesis}[Orbits near the boundary]\label{h2} 
Let $f:R \times R^n \to R^n$ be a $C^\infty$ function. Let $U \subset R \times R^n$ be an open set such that $f$ has only generic bifurcations in $U$.    
Assume also that the set of periodic points
  of $f$ in $U$ is contained in a bounded set in $R^{n+1}$. Let $J$ be
  an interval, and let $Q$ be an open arc in $PO_{nonflip}(f)$ with index
  orientation $\sigma:J \to Q$ such that that $Q$ is contained in $U$,
  though the image of its endpoints may not be. If $[\sigma(\inf J)]$ is defined
  and contained in $V \setminus U$, then it is called an {\bf entry
    orbit} of $U$. Denote the set of all entry orbits of $U$ by $IN$. If
  $[\sigma(\sup J)]$ is defined and is contained in $V \setminus U$,
  then it is called an {\bf exit orbit} of $U$. Denote the set of all exit
  orbits of $U$ by $OUT$. Assume that $IN \cap OUT$ is empty.
\end{hypothesis}

In our case, this hypothesis is satisfied for $F$ for the set $U = \{ [A_0,A_1] \times
E \}$.  In fact, the set of exit orbits for $U$ is empty, and the set
$IN$ of nonflip orbits on the boundary of the region $U$ is entirely
contained in $\{ A_1 \times \mbox{int} (E)\}$.

\begin{theorem}[General cascades theorem]\label{t:main}
  Assume Hypothesis~\ref{h2}. For any odd number $d$, consider the set
  of orbits in $PO_{nonflip}(f)$ with period $2^k d$ for a positive integer $k$.  Let $IN_d$ be
  the set of entry periodic orbits of this type in
  $PO_{nonflip}(f)$. Let $OUT_d$ be the set of exit periodic orbits of
  this type in $PO_{nonflip}(f)$. Assume that $IN_d$ contains $K$
  elements, and $OUT_d$ contains $J$ elements. We allow one but not
  both of $J$ and $K$ to be infinite.
  If $K<J$, then all but $K$ members of $OUT_d$ are contained in
  distinct period-doubling cascades.  Likewise, if $J<K$, then all but
  $J$ members of $IN_d$ are contained in distinct period-doubling cascades.
\end{theorem}

This theorem is proved in~\cite{sander:yorke:p09}. 
From this abstract theorem, we conclude that each nonflip periodic
orbit $P$ for our $F(A_1, \cdot, \cdot)$ is contained in a unique cascade.  We
have already shown that there is a unique nonflip periodic orbit for
$F(A_1,\cdot,\cdot)$ corresponding to each even orbit for the two-shift.

\bigskip
\noindent
{\bf Proof of Step 5: Period of the cascades.}

The only type of generic bifurcations which change the period of an
orbit are period-doubling and period-halving bifurcations. If the
period $k$ of $P$ is odd, then the period is already minimal
in the sense that it is not possible to bifurcate to half the
period. Therefore the cascade through $P$ is period $k$. If the period
$k$ is even, then it is possible that there is a period-halving
bifurcation within the cascade, implying that the period of the
cascade is less than the period $k$ of the orbit $P$. However, the
period always changes by a factor of two, implying that the ratio of
the period of $P$ and the period of the cascade is  a power of
two.

This completes the proof of the theorem.

\end{proof}

\subsection{Reduced smoothness conjecture}\label{s:conjecture}

We end with a conjecture about the non-generic version of our abstract
theorem. If proved, in the context of the current paper, this would
imply that it is possible to extend the results on perturbed H{\'e}non
families to the case of non-generic perturbations. The barrier to
proving Theorem~\ref{t:main} in the non-generic case, is that it is
not possible to control the behavior of eigenvalues other than those
involved in the bifurcations. Such control is needed for the limiting 
arguments to work. Otherwise, a sequence of maps with periodic orbits
of fixed period can limit to a map with a periodic orbit of
smaller period. Lefschetz number arguments, as used by
Franks~\cite{franks:85}, do not apply to Theorem~\ref{t:main},
since the hypotheses use the orbit index, and thus do not include any
assumptions on the Morse index of the flip orbits. The following is a
generalization of the definition of orientation.

\begin{definition}[Generalized orientation]
  Assume that $f:R \times R^n \to R^n$ is a continuous family  such
  that there are only hyperbolic periodic
  orbits on the boundary of a region $U$ in $R \times R^n$. We call an orbit $[q]$ 
   {\bf a generalized entry orbit} at if the Morse index of
  $q$ is even and {\bf a generalized exit orbit} if the Morse
  index is odd.
\end{definition}

Every orbit which is an entry (exit) orbit is a generalized entry
(exit) orbit, but the generalized orientation is also defined for flip
orbits.

\begin{conjecture}[Abstract result reformulated]\label{c:generalized}
  Assume that $f:R \times R^n \to R^n$ is a continuous  family and that $U \in R \times R^n$ is  
such that the periodic orbits in $U$  are contained in a bounded
  region. Assume also that on the boundary of  $U$, there
  are only hyperbolic periodic orbits. Let $IN_d$ and $OUT_d$
  respectively be the generalized entry and exit orbits on the
  boundary of $U$ of period $p=2^m d$, where $d$ is a fixed odd number, and $m$ is any positive integer. Assume that
  the number of orbits in $OUT_d$ and the number of orbits in $IN_d$
  differ.  Let $K$ be the smaller of these two numbers. (We allow one
  but not both of these numbers to be infinite.)  Then there is a
  cascade through all but possibly $K$ of the orbits in the larger of the sets
  $IN_d$ and $OUT_d$.
\end{conjecture}

\section{Acknowledgements}

E.S. was partially supported by NSF Grant DMS-0639300
and NIH Grant R01-MH79502. J.A.Y. was partially supported by NSF Grant DMS-0616585 and NIH Grant
R01-HG0294501.

\begin{flushleft}
%
% References
%
\addcontentsline{toc}{subsection}{References}
\footnotesize
%
%   Add in the bbl-file the command \parskip 0pt.
%

%

{\bf AMS Subject Classification: 37.}\\[2ex]

% Write more than one author separately if they have different 
% affiliations, otherwise write the names on the same line, separated 
% by commas.
%
E.~Sander
Department of Mathematical Sciences,
George Mason University,
4400 University Dr.,
Fairfax, VA, 22030, USA. 
E-mail: \texttt{sander@math.gmu.edu}

J.A.~Yorke
Department of Mathematics, IPST, and Physics Department,
University of Maryland,
College Park, MD 20742, USA. 
E-mail: \texttt{yorke@umd.edu}

\end{flushleft}

\end{document}